# A NEW MULTIPLE TESTING METHOD IN THE DEPENDENT CASE[1]

By Arthur Cohen, Harold B. Sackrowitz and Minya Xu

*Rutgers University*

The most popular multiple testing procedures are stepwise procedures based on $P$-values for individual test statistics. Included among these are the false discovery rate (FDR) controlling procedures of Benjamini–Hochberg [*J. Roy. Statist. Soc. Ser. B* **57** (1995) 289–300] and their offsprings. Even for models that entail dependent data, $P$-values based on marginal distributions are used. Unlike such methods, the new method takes dependency into account at all stages. Furthermore, the $P$-value procedures often lack an intuitive convexity property, which is needed for admissibility. Still further, the new methodology is computationally feasible. If the number of tests is large and the proportion of true alternatives is less than say 25 percent, simulations demonstrate a clear preference for the new methodology. Applications are detailed for models such as testing treatments against control (or any intraclass correlation model), testing for change points and testing means when correlation is successive.

**1. Introduction.** The need for multiple testing procedures (MTPs) has been given great impetus by diverse fields of application such as microarrays, astronomy, mutual fund evaluations, proteomics, disclosure risk, cytometry, imaging and others. Traditional methods to deal with multiple testing when the number of tests is large are deemed too conservative (i.e., they do not detect significant effects often enough). New approaches to multiple testing have arisen. Many of the new approaches are classified as stepwise procedures, such as step-up or step-down in contrast to single step procedures [see Hochberg and Tamhane (1987) and also Dudoit, Shaffer and Boldrick (2003), where 18 procedures are listed as single step, step-up or step-down].

Received July 2007; revised April 2008.
[1]Supported by NSF Grant DMS-04-57248 and NSA Grant H-98230-06-0076.

Among the more popular procedures is the Benjamini–Hochberg (1995) false discovery rate (FDR) controlling procedure. Many offsprings have followed [see, e.g., Efron et al. (2001), Storey and Tibshirani (2003), Sarkar (2002), Benjamini and Yekutieli (2001) and Cai and Sarkar (2006), just to mention a few]. Typically, the stepwise procedures deal with $P$-values determined from marginal distributions [see, e.g., Dudoit and van der Laan (2008), Chapter 3]. Even when the model entails random vectors with correlated variates, $P$-values from marginal distributions, ignoring correlations, are the basis of the procedures.

Many multiple testing procedures are designed to control some error rate such as the familywise error rate FWER (weak and strong), $k$-FWER [see Lehmann and Romano (2005)] and FDR. However, many researchers study the multiple testing problem as a finite action decision problem with a variety of loss functions [see, e.g., Lehmann (1957), Genovese and Wasserman (2002), Ishwaran and Rao (2003) and Muller et al. (2004)]. In these studies, the merits of the procedures are evaluated and compared by their risk functions. The risk function approach does not always necessitate the need to control a particular kind of error rate and can sometimes lead to procedures whose overall performance is preferred or even strongly preferred to an error controlling procedure. Whereas FDR control is appropriate for some situations where the number of tests is large, there are many situations where one would prefer a procedure whose expected number of both type I errors and type II errors are smaller. Dudoit and van der Laan (2008) study expected values of functions of numbers of type I and type II errors.

In a series of papers [Cohen and Sackrowitz (2005, 2007, 2008) and Cohen, Kolassa and Sackrowitz (2007)] demonstrated that, given a typical step-up or step-down procedure, there exist other procedures whose expected numbers of type I and type II errors are smaller. In fact, in Cohen and Sackrowitz (2008), for multivariate normal models when correlation is nonzero for two-sided alternatives, there exist procedures whose individual tests have smaller expected type I and type II errors.

In this paper, we assume $\mathbf{X}$ is an $M \times 1$ vector that is multivariate normal with mean vector $\boldsymbol{\mu}$ and covariance matrix $\Gamma = \sigma^2\Sigma$. The matrix $\Sigma$ is a known positive, definite nondiagonal matrix. The parameter $\sigma^2$ is either known or unknown. In the latter case an estimator of $\sigma^2$, which is a scaled chi-square variable, is available and this variable is independent of $\mathbf{X}$. This is a classical linear model assumption. $\Sigma$ is known since it is a function of the design matrix. We will demonstrate the new methodology in two important subclasses of this model. The first is the intraclass covariance matrix model, which characterizes the popular situation in which the variables are exchangeable. This model includes the problem of testing several treatments against a control. The second application is to the successive correlation covariance matrix, which includes change point problems. We test two sided alternatives (i.e., $H_i : \mu_i = 0$ vs. $K_i : \mu_i \neq 0$, $i = 1, \ldots, M$). We also test one sided alternatives (i.e., $H_i^* : \mu_i \leq 0$ vs. $K_i^* : \mu_i > 0$, $i = 1, \ldots, M$ or $H_i : \mu = 0$ vs. $K_i^* : \mu_i > 0$).

The goal of this paper is to develop good MTPs in the case of correlated variables. To begin with, we realize that every MTP induces individual tests, $\phi_i$, for the individual hypothesis testing problems $H_i$ vs. $K_i$. The behavior of these tests should be of fundamental concern. However, the stepwise construction of most MTPs often makes it difficult to describe and study the individual tests.

In particular, suppose an individual test induced by an MTP is inadmissible for the standard hypothesis testing loss. That is, for that individual hypothesis testing problem, a test exists whose size is no greater than the stepwise procedure test and whose power is no less with some strict inequality. It would then follow that the overall procedure would be inadmissible whenever the risk function is a monotone function of the expected numbers of type I and type II errors.

As a first step, we find a convexity property that is necessary and sufficient for admissibility of the individual tests. In Cohen and Sackrowitz (2008), it has been shown that most popular stepwise procedures do not possess the convexity property when there is correlation in the two-sided alternative case. Next, we construct a step-down type MTP whose individual tests do have the required convexity property. As is typical in problems where no single optimal procedure exists, the selection of a procedure is somewhat subjective. In evaluating procedures, we focus mainly on the expected number of type I and type II errors that the procedures make.

The new stepwise testing method proposed is based on the maximum of adaptively formed residuals. The method is called maximum residual down (MRD). The MRD method has several advantages over the stepwise methods that are currently recommended in the literature:

(1) The main justification for MRDs is the fact that MRD tests take into account the correlation among the $M$ variates. Thus, MRD utilizes information oftentimes not used by the current $P$-value methods. This property of the MRD procedure is the likely explanation for the apparent improved overall performance of MRD when compared to the $P$-value methods based on marginal distributions.

(2) MRD procedures have an intuitive and desirable convexity property required for admissibility. Whereas admissibility is not in itself a compelling property, inadmissibility can be a serious shortcoming.

(3) For the treatment vs. control and change point models for large and relevant portions of the parameter space, simulations demonstrate that the MRD method makes substantially fewer mistakes than the popular FDR controlling procedures. In particular, if the proportion of true alternatives is less than 25 percent of the total number of tests, then the simulations are somewhat convincing that this method is quite good.

(4) The MRD method is applicable in all cases where $\Sigma$ is known.

For arbitrary $\Sigma$ and $M$ extremely large, the procedure essentially requires inversion of a larger size matrix. This could be computationally difficult. The level of difficulty depends on the structure of $\Sigma$. In the two popular models, we consider $\Sigma$ can easily be inverted regardless of how large $M$ is. The first model is intraclass. The concept of intraclass covariance matrix was introduced by Rao (1945). Subsequently it has been discussed in articles in behaviorial genetics and statistics [see, e.g., Carey (2005) and Krishnaiah and Pathak (1967)]. Such a model is appropriate whenever the components of $\mathbf{X}$ have a multivariate distribution that is exchangeable. In particular, all variances are equal and all covariances are equal. The intraclass correlation matrix is appropriate for the model of testing each of $M-1$ treatments against a control. MRD is readily applicable here, since inversion of the appropriate covariance matrix is easily facilitated. The second model is successive correlation. This model has a constant nonzero correlation coefficient between adjacent pairs of variables. All other correlations are zero [see Krishnaiah and Pathak (1967)]. This model presents no computational issues even if $M$ is extremely large. A special case of this model is the change point problem [see Chen and Gupta (2000)]. We will see in Section 6 that the MTP method discussed in Chen and Gupta (2000) is based on many collections of pooled means. This is precisely the set of statistics given by the MRD method applied to this very special case. In a sense, this validates our very general approach, even though our method uses the statistics differently than in Chen and Gupta (2000).

A seemingly logical step-down method that would take correlations into account is to successively perform likelihood ratio tests (LRT) of global hypotheses. That is, one could employ the closure method [see Marcus, Peritz and Gabriel (1976)] using an LRT, at step one, for $\boldsymbol{\mu}=\mathbf{0}$ vs. $\boldsymbol{\mu}\neq 0$. If the global test rejects, then eliminate the variate corresponding to $\max_{1\le i\le M}|X_i|$. One continues in a step-down fashion in determining the LRT-based MTP. Call this procedure LRSD.

When $\Sigma$ is intraclass for two-sided alternatives, LRSD is admissible for any monotone collection of critical constants only when $M=2$ or $M=3$. For $M\ge 4$, counterexamples abound. That is, there are many critical constants for which LRSD is inadmissible. Furthermore, critical constants are found for $M\ge 5$, which relate to constants that are likely to be used. This inadmissibility of LRSD is what prompted and led us to MRD.

For one-sided alternatives when $\Sigma$ is intraclass, LRSD is admissible even in cases where the common variance $\sigma^2$ is unknown (provided replications of the observations are taken). In this instance, LRSD can be a competitor to MRD, and this is reflected in the simulation study in Section 7.

We note there that $M$ is taken to be 100. Large values of $M$ entail computational problems for LRSD, since under the alternative the parameter space is constrained, and the software needed to carry out the tests for $M$ much greater than 100 is very time consuming. For the treatments vs. control model, one might think that the $P$-value based step-down procedure, based on the analogue of Dunnett's one sided tests of global hypotheses, might also be a competitor [see Westfall and Young (1993), Section 3.2.1]. In this instance, dependency is taken into account when determining critical values. Nevertheless, it does not take correlation into account in the test statistics, and, overall, the procedure does not fare well in the simulation study.

As previously mentioned, another problem of interest is testing for change points in a sequence of $M+1$ independent normal trials. Sometimes, it is assumed that the means are nondecreasing, in which case the model is referred to as a simple order model. One seeks to determine whether a change in mean has occurred at particular time points. The alternative at each time point is either two-sided or one-sided. In either case, the LRSD step-down method is mostly inadmissible, while the MRD method is admissible.

Returning to the general case, we remark that if $\Sigma$ is unknown but replications are available, an estimator of $\Sigma$ can replace it in the MRD method. We cannot claim the optimality properties, but, nevertheless, the method is viable. For large numbers of replications, even the normality assumption may not be crucial.

In the next section, we describe the MRD method. In Section 3, we prove that the MRD method is admissible for the vector risk where each component of the vector is the testing risk for an individual test. Admissibility depends on whether each individual test function has an intuitive convexity property. Section 4 is concerned with the LRT based step-down procedure (LRSD). Here, there are both admissibility and inadmissibility results of interest. Section 5 contains a geometric connection between the MRD and LRSD methods, some other interesting interpretations and some figures related to the geometric interpretation. Results concerned with testing several treatments vs. control, the change point problem and the successive correlation model are given in Section 6. Simulations and analyses are given in Section 7. Most proofs appear in the Appendix.

## 2. MRD method.

Assume $\mathbf{X} = (X_1, \ldots, X_M)'$ is distributed according to a multivariate normal distribution with mean vector $\boldsymbol{\mu}$ and covariance matrix $\sigma^2\Sigma = \sigma^2(\sigma_{ij})$. The matrix $\Sigma$ is assumed known, and, for now, we take $\sigma^2$ to be known and, without loss of generality, let $\sigma^2 = 1$. The two-sided multiple testing problem is test

$$H_i : \mu_i = 0 \quad \text{vs.} \quad K_i : \mu_i \neq 0, \qquad i = 1, \ldots, M. \tag{2.1}$$

We will also consider one-sided alternative problems

$$H_i : \mu_i = 0 \quad \text{vs.} \quad K_i^* : \mu_i > 0 \tag{2.2}$$

and

$$H_i^* : \mu_i \le 0 \quad \text{vs.} \quad K_i^* : \mu_i > 0. \tag{2.3}$$

For now, we focus on the two-sided case (2.1).

By way of notation, $\mathbf{X}^{(i_1,i_2,\ldots,i_r)}$ is the $(M-r)$ vector consisting of the components of $\mathbf{X}$ with $X_{i_1},\ldots,X_{i_r}$ left out. $\Sigma_{(i_1,\ldots,i_r)}$ is the $(M-r)\times(M-r)$ covariance matrix of $\mathbf{X}^{(i_1,\ldots,i_r)}$. $\boldsymbol{\sigma}_{(j)}^{(i_1,\ldots,i_{m-1})}$ is the $(M-m)\times 1$ vector of covariances between $X_j$ and all variables except $X_{i_1},\ldots,X_{i_{(m-1)}}$, and $X_j$.

$$\sigma_{(j\cdot i_1,\ldots,i_{(m-1)})} = \sigma_{jj} - \boldsymbol{\sigma}_{(j)}^{(i_1,\ldots,i_{(m-1)})'} \Sigma^{-1}_{(i_1,\ldots,i_{(m-1)},j)} \boldsymbol{\sigma}_{(j)}^{(i_1,\ldots,i_{(m-1)})}$$

is the conditional variance of $X_j$, given all variables except $X_{i_1},\ldots,X_{i_{(m-1)}}, X_j$.

Now, define

$$\begin{aligned} &U_{mj}^{(i_1,\ldots,i_{m-1})}(\mathbf{X}) \\ &\quad = (X_j - \boldsymbol{\sigma}_{(j)}^{(i_1,\ldots,i_{(m-1)})'} \Sigma^{-1}_{(i_1,\ldots,i_{(m-1)},j)} \mathbf{X}^{(i_1,\ldots,i_{(m-1)},j)})/\sigma^{1/2}_{(j\cdot i_1,\ldots,i_{(m-1)})} \end{aligned} \tag{2.4}$$

for $m = 1,\ldots,M$.

The $m$ subscript represents the stage of the MRD procedure. Note that

$$U_{mj}^{(i_1,\ldots,i_{m-1})} = (X_j - E_0\{X_j|\mathbf{X}^{(i_1,\ldots,i_{(m-1)},j)}\})/\sqrt{\operatorname{Var}(X_j|\mathbf{X}^{(i_1,\ldots,i_{(m-1)},j)})},$$

where $E_0$ is taken under $\boldsymbol{\mu} = 0$.

We now describe a general class of stepwise down procedures, given a set of $M2^{M-1}$ functions $U_{mj}(\mathbf{x})$. At most, $M(M+1)/2$ of these needs to be calculated to carry out the procedure. The $m$ index ranges from $1,2,\ldots,M$ and represents the $m$th stage. At stage $m$ there are $M-m+1$ functions.

Let $C_1 > C_2 > \cdots > C_M > 0$ be a given set of constants. At stage 1, consider $U_{1j}(\mathbf{x})$, $j \in \{1,\ldots,M\}$. Let $j_1 = j_1(\mathbf{x})$ be such that $U_{1j_1}(\mathbf{x}) = \max_j |U_{1j}(\mathbf{x})|$. If $U_{1j_1}(\mathbf{x}) < C_1$, stop and accept all $H_i$. Otherwise, reject $H_{j_1}$ and continue to stage 2.

At stage 2, consider $M-1$ functions $U_{2j}^{(j_1)}(\mathbf{x})$, $j \in \{1,\ldots,M\}\setminus\{j_1\}$. Note that $U_{2j}^{(j_1)}(\mathbf{x})$ just depends on $\mathbf{x}^{(j_1)}$. Let $j_2 = j_2(\mathbf{x}^{(j_1)})$ be such that $U_{2j_2} = \max_j |U_{2j}^{(j_1)}|$, $j \in \{1,\ldots,M\}\setminus\{j_1\}$. If $U_{2j_2} < C_2$, stop and accept all remaining null hypotheses. Otherwise, reject $H_{j_2}$ and continue to stage 3.

In general, at stage $m$, $m = 1,\ldots,M$, consider $M-m+1$ functions $U_{mj}^{(j_1,\ldots,j_{(m-1)})}(\mathbf{x})$, $j \in \{1,\ldots,M\}\setminus\{j_1,\ldots,j_{(m-1)}\}$. Note that $U_{mj}^{(j_1,\ldots,j_{(m-1)})}$ depends on $\mathbf{x}^{(j_1,\ldots,j_{(m-1)})}$. Let $j_m = j_m(\mathbf{x}^{(j_1,\ldots,j_{(m-1)})})$ be such that

$U_{mj_m} = \max_j |U_{m_j}^{(j_1,\dots,j_{(m-1)})}|$, $j \in \{1,\dots,M\} \setminus \{j_1,\dots,j_{(m-1)}\}$. If $U_{mj_m} < C_m$, stop and accept all remaining null hypotheses. Otherwise, reject $H_{j_m}$ and continue to stage $m+1$ (unless $m = M$, in which case, stop).

The above MTP determines test functions for each individual testing problem. Let $\phi_U(\mathbf{x})$ denote the test function for testing $H_1 : \mu_1 = 0$ vs. $K_1 : \mu_1 \neq 0$.

Note that, at the beginning of Section 7, we offer a discussion regarding the choice of $C_1,\dots,C_m$.

**3. Admissibility of MRD.** We will demonstrate that for each individual testing problem that the MTP based on the MRD method is admissible. Without loss of generality, we focus on $H_1$ vs. $K_1$ and start with the case $\sigma$ known ($\sigma^2 = 1$). Our plan is to use a result of Matthes and Truax (1967) which offers a necessary and sufficient condition for admissibility of a test of $H_1$ vs. $K_1$ when the joint distribution of $\mathbf{X}$ is an exponential family. We will state this result for our model as Lemma 3.1. We next demonstrate, in Lemma 3.2, that the $U_{mj}$ function given in (2.4) has certain monotonicity properties. These monotonicity properties will enable us to prove, in Lemma 3.3, that the individual test functions for $H_i$ vs. $K_i$ have a convexity property that is necessary and sufficient for admissibility. Theorem 3.1 summarizes and states the admissibility of the MRD procedure.

Now, we express the density of $\mathbf{X}$ as

$$f_{\mathbf{X}}(\mathbf{x}|\boldsymbol{\mu}) = (1/(2\pi)^{M/2}|\Sigma|^{1/2}) \exp -\tfrac{1}{2}(\mathbf{x}-\boldsymbol{\mu})'\Sigma^{-1}(\mathbf{x}-\boldsymbol{\mu}), \tag{3.1}$$

which, in exponential family form, is

$$f_{\mathbf{X}}(\mathbf{x}|\boldsymbol{\mu}) = h(\mathbf{x})\beta(\boldsymbol{\mu}) \exp \mathbf{x}'\Sigma^{-1}\boldsymbol{\mu}. \tag{3.2}$$

Next, let $\mathbf{Y} = \Sigma^{-1}\mathbf{X}$ so that

$$f_{\mathbf{Y}}(\mathbf{y}|\boldsymbol{\mu}) = h^*(\mathbf{y})\beta(\boldsymbol{\mu}) \exp \sum_{i=1}^{M} y_i\mu_i. \tag{3.3}$$

LEMMA 3.1. *A necessary and sufficient condition for a test $\phi(\mathbf{y})$ of $H_1 : \mu_1 = 0$ vs. $K_1 : \mu_1 \neq 0$ to be admissible is that, for almost every fixed $y_2,\dots,y_M$, the acceptance region of the test is an interval in $y_1$.*

PROOF. See Matthes and Truax (1967). □

Note that, to study the test function $\phi(\mathbf{y}) = \phi_U(\mathbf{x})$ as $y_1$ varies and $(y_2,\dots,y_M)$ remain fixed, we can consider sample points $\mathbf{x} + r\mathbf{g}$ where $\mathbf{g}$ is the first column of $\Sigma$ and $r$ varies. This is true, since $\mathbf{y}$ is a function of $\mathbf{x}$, and so $\mathbf{y}$ evaluated at $(\mathbf{x} + r\mathbf{g})$ is $\Sigma^{-1}(\mathbf{x} + r\mathbf{g}) = \mathbf{y} + (r, 0, \dots, 0)' = (y_1 + r, y_2, \dots, y_M)'$.

From here on it will be convenient to express the functions $U_{mj}^{(j_1,\dots,j_{m-1})}(\mathbf{x})$ of (2.4) simply as $U_{mj}(\mathbf{x})$. No confusion should ensue.

LEMMA 3.2. *The functions $U_{mj}(\mathbf{x})$ given in (2.4) have the following properties.*

*As a function of $r$,*

$$U_{m1}(\mathbf{x}+r\mathbf{g}) = U_{m1}(\mathbf{x}) + r. \tag{3.4}$$

*For $m = 1,\dots,M$; $j \in \{2,\dots,M\} \setminus \{j_1,\dots,j_{m-1}\}$, $j_1 \neq 1,\dots,j_{m-1} \neq 1$,*

$$U_{mj}(\mathbf{x}+r\mathbf{g}) = U_{mj}(\mathbf{x}). \tag{3.5}$$

PROOF. See Appendix. □

LEMMA 3.3. *Suppose that, for some $\mathbf{x}^*$ and $r_0 > 0$, $\phi_U(\mathbf{x}^*) = 0$ and $\phi_U(\mathbf{x}^* + r_0\mathbf{g}) = 1$. Then, $\phi_U(\mathbf{x}^* + r\mathbf{g}) = 1$ for all $r > r_0$.*

PROOF. See Appendix. □

Note that Lemma 3.3 implies that the acceptance region in $y_1$, for fixed $y_2,\dots,y_m$ is an interval.

THEOREM 3.1. *For the two sided case, the MRD procedure based on $\{U_{mj}\}$ is admissible.*

PROOF. Admissible means that each individual test for each hypothesis testing problem is admissible. Without loss of generality, we show admissibility of $\phi_U(\mathbf{x})$ for $H_1$ vs. $K_1$. Proof that the other tests are admissible for the other hypotheses would be done the same way. That $\phi_U(\mathbf{x})$ is admissible for $H_1$ vs. $K_1$ follows readily from Lemmas 3.1 and 3.3. □

For the case where $\sigma^2$ is unknown, we assume that we have available an unbiased estimator $s^2$ of $\sigma^2$ with the property that $\nu s^2/\sigma^2$ is a $\mathcal{X}_\nu^2$ variable that is independent of $\mathbf{X}$. In this situation, we write the joint density of $(\mathbf{X}, s^2)$

$$\begin{aligned} &f_{\mathbf{X},s^2}(\mathbf{x}, s^2|\mu, \sigma^2) \\ &\quad = (\nu(\nu s^2)^{\nu/2-1}/(2\pi)^{M/2} \cdot 2^{\nu/2}\Gamma(\nu/2)(\sigma^2)^{(M+\nu)/2}|\Sigma|^{1/2}) \\ &\qquad \times \exp(-1/2\sigma^2)\{(\mathbf{x}-\boldsymbol{\mu})'\Sigma^{-1}(\mathbf{x}-\boldsymbol{\mu}) + \nu s^2\} \\ &\quad = H(\mathbf{x}, T)B(\mu, \sigma^2)\exp\{\mathbf{x}'\Sigma^{-1}\boldsymbol{\mu}/\sigma^2 - (1/2\sigma^2)T\}, \end{aligned} \tag{3.6}$$

where $T = \nu s^2 + \mathbf{x}'\Sigma^{-1}\mathbf{x}$. Note that the change from $(\mathbf{x}, s^2)$ to $(\mathbf{x}, T)$ limits the values of $\mathbf{x}$ in the sample space to those for which $\mathbf{x}'\Sigma^{-1}\mathbf{x} \le T$.

The MRD method now utilizes statistics $U_{mj}(\mathbf{X})/T^{1/2}$. All lemmas and Theorem 3.1 hold, with $T$ as well as $y_2, \ldots, y_M$ fixed, where, once again, $\mathbf{y} = \Sigma^{-1}\mathbf{x}$.

For one-sided alternatives specified in (2.2) and (2.3), the MRD method simply uses $U_{mj}$ in place of $|U_{mj}|$. The result of Theorem 3.1 can be proved similarly.

**4. Likelihood ratio step-down method (LRSD).** Assume that $\mathbf{X}$ is distributed as multivariate normal with unknown mean vector $\boldsymbol{\mu}$ and known intraclass covariance matrix $\Sigma$. Without loss of generality, we take the diagonal elements of $\Sigma$ to be 1 and the off diagonal elements to be $\rho$. The LRSD method is to test by the LRT, at stage 1, the global hypothesis $H_{1G}: \boldsymbol{\mu} = \mathbf{0}$ vs. $K_{1G}: \boldsymbol{\mu} \neq 0$. If $H_{1G}$ is not rejected, then stop and accept all $H_i$, $i = 1, \ldots, M$. If $H_{1G}$ is rejected, then reject $H_{j_1}$, where $j_1$ is the index for which $|X_{j_1}| = \max_{1 \le j \le M} |X_j|$, and continue to stage 2. At stage 1 use the critical value $C_1$. At stage 2, test, by the LRT the global hypothesis, $H_{2G}: \boldsymbol{\mu}^{(j_1)} = \mathbf{0}$ vs. $K_{2G}: \boldsymbol{\mu}^{(j_1)} \neq 0$, where $\boldsymbol{\mu}^{(j_1)}$ is the $(M-1) \times 1$ vector of means that are the same as $\boldsymbol{\mu}$ save $\mu_{j_1}$ is left out. Use the critical value $C_2 < C_1$. Proceed as in stage 1. At stage $m$, test by the LRT $H_{mG}: \boldsymbol{\mu}^{(j_1, \ldots, j_{m-1})} = 0$ vs. $K_{mG}: \boldsymbol{\mu}^{(j_1, \ldots, j_{m-1})} \neq \mathbf{0}$ and so on. At stage $m$, use the critical value $C_m < C_{m-1}$.

We will demonstrate that the LRSD is admissible for $M = 2$ and $M = 3$. For $M \ge 4$ there exist counterexamples for certain collections of critical values and certain values of $\rho$. We offer a counterexample when $M = 4$, and, when $M = 5$, we demonstrate inadmissibility for a large class of practical critical values for logical values of $\rho$. In fact, for large $M$ using $\chi^2$ critical values, it turns out that for most $\rho$ values ($\rho \neq 0$) counterexamples demonstrate that LRSD is inadmissible.

On the other hand, should the alternatives for the individual hypotheses be the one-sided alternatives given in (2.2), then the LRSD is admissible.

When the alternative is two-sided, the results of Section 3 imply that admissibility of a test for an individual hypothesis testing problem (say $H_1$ vs. $K_1$, without loss of generality) is determined by whether the conditional acceptance region in $y_1$ given $(y_2, \ldots, y_M)$ is an interval. (Recall $\mathbf{y} = \Sigma^{-1}\mathbf{x}$.) When the alternative is one-sided, the conditional acceptance region is a left sided half line.

Focusing first on the two-sided alternative case, we note that the LRT for $H_{1G}$ vs. $K_{1G}$ is to reject if

$$\mathbf{x}'\Sigma^{-1}\mathbf{x} \ge C_1, \tag{4.1}$$

where

$$\Sigma^{-1} = (1/(1-\rho))\{I - G(\mathbf{11}')\} \tag{4.2}$$

and $G = \rho/(1+(M-1)\rho)$. As such,

$$\Sigma^{-1} = [1/(1-\rho)(1+(M-1)\rho)] \times \begin{pmatrix} 1+(M-2)\rho & -\rho & \cdots & \cdots & \cdots & -\rho \\ -\rho & (1+(M-2)\rho) & & -\rho & \cdots & -\rho \\ -\rho & -\rho & \cdots & -\rho & & (1+(M-2)\rho) \end{pmatrix}. \tag{4.3}$$

Again, let $\mathbf{g}$ be the first column of $\Sigma$.

Our first result in this section follows.

THEOREM 4.1. *For the two-sided alternative case, LRSD is admissible for $M = 2$ and $M = 3$.*

PROOF. See Appendix. □

For $M = 4$, we exhibit a set of critical values for which LRSD is inadmissible. To do so, we find a sample point $\mathbf{x}^*$ at which $H_1$ is rejected and for which $H_1$ is accepted at $\mathbf{x}^* + \gamma\mathbf{g}$. In fact, let $\mathbf{x}^* = (a, -a-\Delta, b, -b-\epsilon)'$ for $b > a+\Delta > a > 0$ and $\epsilon > 0$. Thus, using (4.1) at stage 1, choose $C_1$ so that $\mathbf{x}^*\Sigma^{-1}\mathbf{x}^* > C_1$ so that $H_4$ is rejected and variable $x_4^*$ is eliminated at stage 2. At stage 2, we calculate

$$\begin{aligned} &\mathbf{x}^{*(4)'}\Sigma_{(4)}^{-1}\mathbf{x}^{*(4)} \\ &\quad = [1/(1+\rho-2\rho^2)] \\ &\qquad \times \{(1+\rho)b^2 + 2a^2(1+2\rho) + 2\Delta[a + 2a\rho + \rho b + (1+\rho)\Delta/2]\}. \end{aligned} \tag{4.4}$$

We set $\mathbf{x}^{*(4)'}\Sigma_{(4)}^{-1}\mathbf{x}^{*(4)} = C_2$. At stage 3, $H_2$ is rejected, and, at stage 4, $H_1$ is rejected. Now, if $\rho > 0$, let $\gamma = \epsilon/\rho$ and note that $(\mathbf{x}^* + \gamma\mathbf{g})'\Sigma^{-1}(\mathbf{x}^* + \gamma\mathbf{g}) > C_1$. This time, however, $H_3$ is rejected at stage 1. At stage 2, we calculate, for $\gamma = \epsilon/\rho$,

$$(\mathbf{x}^{*(3)} + \gamma\mathbf{g}^{(3)})'\Sigma_{(3)}^{-1}(\mathbf{x}^{*(3)} + \gamma\mathbf{g}^{(3)}). \tag{4.5}$$

We note that (4.4) minus (4.5) is

$$\begin{aligned} &1/(1+\rho-2\rho^2)\{4\Delta b\rho - \varepsilon^2[\rho - 1 + 1/\rho + 1/\rho^2] \\ &\qquad - \varepsilon[2a/\rho + 2a - 4a\rho - 2\rho\Delta + 2(1+\rho)b]\}. \end{aligned} \tag{4.6}$$

There are many choices of $a, b, \Delta, \epsilon, \rho, \gamma$ for which (4.6) is positive (e.g., $a = 2, b = 4, \Delta = 1, \epsilon = 0.1, \rho = 0.5, \gamma = 0.2$). The fact that (4.6) $> 0$ implies that at $\mathbf{x}^* + \gamma\mathbf{g}$ the overall procedure rejects $H_3$ and accepts $H_1, H_2$ and $H_4$. Note that, since $x_1^* > 0$, $\mathbf{x}^* - x_1\mathbf{g}$ is an accept point. Now, if $H_1$ is rejected for $\mathbf{x} = \mathbf{x}^*$ but accepted for $\mathbf{x}^* + \gamma\mathbf{g}$, it is implied that the test for $H_1$ is inadmissible.

For $M = 5$, it can be shown that if the critical values correspond to critical values of chi-square with $m$ degrees of freedom, $m = 1, 2, 3, 4, 5$ at level, say 0.05, then, for most values of $\rho$, LRSD is also inadmissible. The same is true for any $M \geq 5$.

Next, for the intraclass model, we consider testing one-sided alternatives (i.e., we test $H_i : \mu_i = 0$ vs. $K_i^* : \mu_i > 0$). The LRSD method in this case is the same as in the two-sided alternative case, except that $|X_{j_1}|$ is replaced by $X_{j_1} = \max(X_1, \ldots, X_M)$ and similarly at subsequent stages. For this setup we have the following theorem.

THEOREM 4.2. *For the one-sided alternative case, LRSD is admissible.*

PROOF. See Appendix. □

The final result of this section deals with the intraclass model when the covariance matrix is of the form $\sigma^2\Sigma$, where $\sigma^2$ is unknown and $\Sigma$ is as before. This time, however, a random sample $\mathbf{X}_\alpha$, $\alpha = 1, \ldots, n$, is taken from a normal distribution with mean vector $\boldsymbol{\mu}$ and covariance matrix $\sigma^2\Sigma$. The alternative hypotheses are one-sided, and the global likelihood ratio test is based on $\bar{\mathbf{X}} = \sum_{\alpha=1}^n \mathbf{X}_\alpha / n$ and $T = \sum_{\alpha=1}^n \mathbf{X}'_\alpha \Sigma^{-1} \mathbf{X}_\alpha$. Using the fact that $\bar{\mathbf{X}}, T$ have an exponential family distribution and arguments similar to those used previously, it can be shown that the LRSD procedure is admissible in this case as well.

We remark that this model ensues for the problem of testing $M$ treatments against a control when it is assumed that the mean for each treatment is greater than or equal to the mean for the control. More details are given in Section 6.

**5. Geometric and other interpretations.** LRSD compared statistics of the form $\mathbf{x}'\Sigma^{-1}\mathbf{x}$ to critical values in order to test global hypotheses at each stage of the process. The overall acceptance region for the global testing problem is therefore an ellipsoid. The individual statistics $\pm U_{mj}$ given in (2.4) determine the MRD method. These statistics represent pairs of supporting hyperplanes to the ellipsoids determining acceptance regions of the global hypotheses at stage $m$ [see Scheffé (1959), page 69]. The particular hyperplanes are tangent to the ellipsoids at sample points on the ellipsoid for which all but one coordinates are zero. If the probability, under a global null hypothesis of a mean vector, is zero, specified at, say $\gamma$, then the probability of the ellipsoid is $\gamma$. The acceptance set determined by the supporting hyperplanes would be larger than $\gamma$. However, should one desire this set to have probability $\gamma$, then the hyperplanes would support a smaller ellipse. Figures 1 and 2 depict such sets for the stage when there are only two means left to test.

Note that, by comparing the $\max_j |U_{mj}|$ to a critical constant, one is determining an acceptance region for the global hypothesis at stage $m$ by using the union-intersection procedure [see, e.g., Casella and Berger (2002), page 380].

The statistics $U_{mj}$ appear in the identity in Anderson (1984), Exercise 54, Chapter 2. Thus, one can express the MRD method alternatively in terms of $\mathbf{x}^{(i)\prime}\Sigma_{(i)}^{-1}\mathbf{x}^{(i)}$.

The statistics $U_{mj}$ are also the focal point in determining change points in the methodology offered by Vostrikova (1981).

Although MRD uses $U_{mj}$ as does Vostrikova (1981), the methodologies are different. We discuss this further in Section 6.

Remark 5.1. It is interesting to note that the MRD method is not $P$-value monotone in the sense of Hommel and Bernhard (1999). That is, an MTP is monotone with respect to $P$-values if $P_i \leq P_j$ and $H_i$ is not rejected, then $H_j$ cannot be rejected. As indicated in that reference, $P$-value monotonicity is not always desirable.

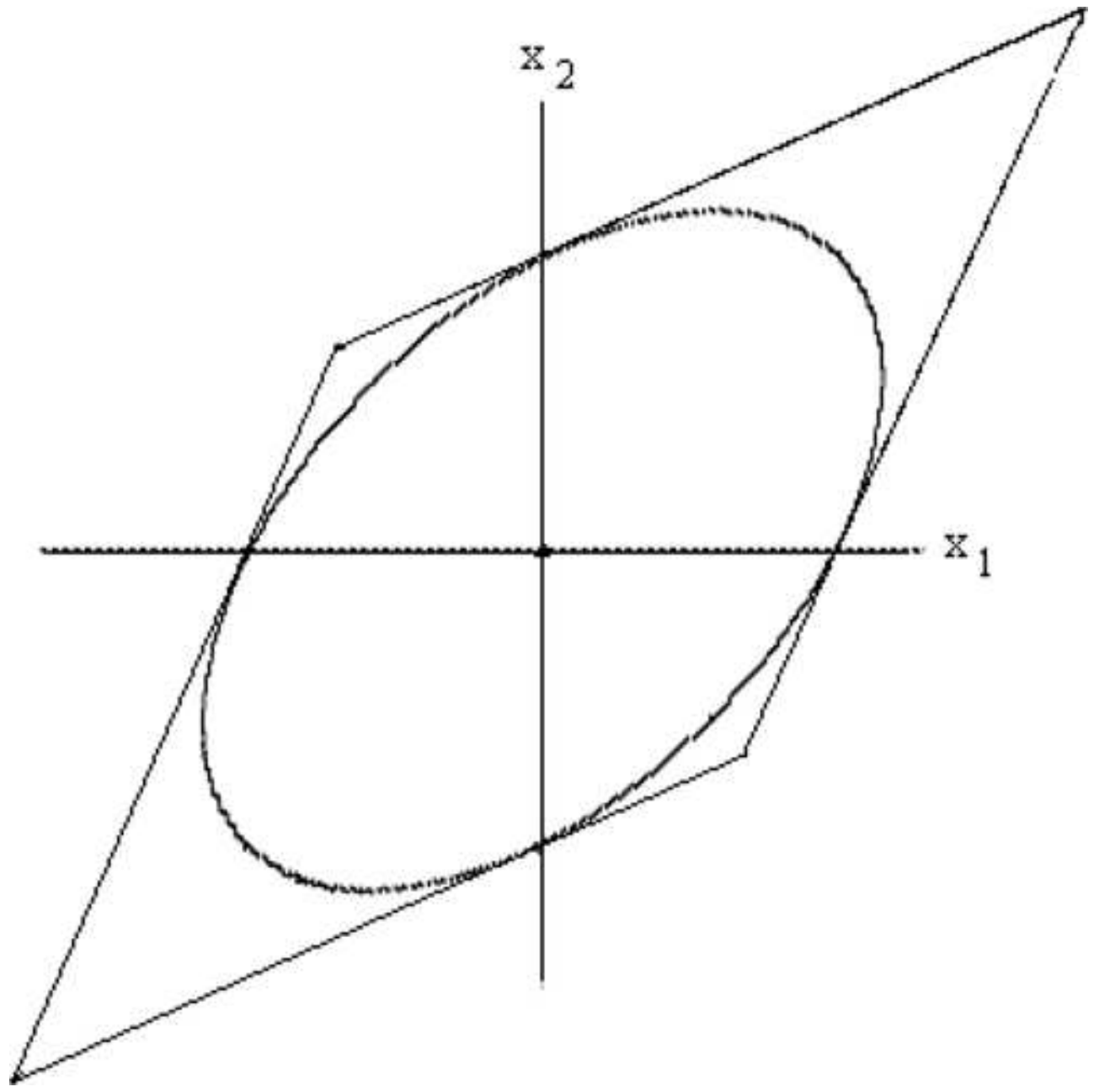

Fig. 1. *LRSD ellipse with supporting hyperplanes in two dimensions.*

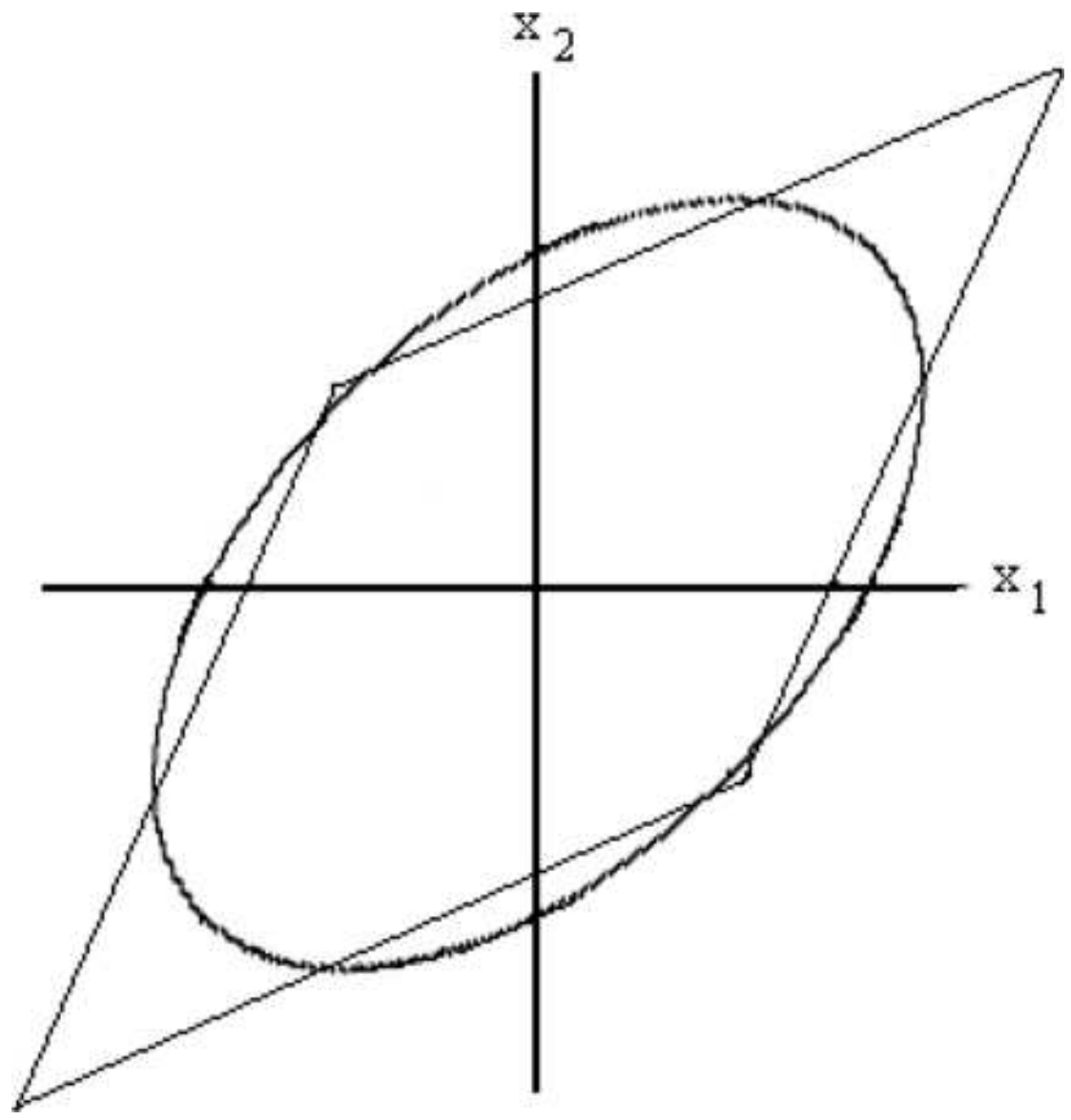

FIG. 2. *LRSD ellipse with supporting hyperplanes shrunk to match size.*

REMARK 5.2. Use of $U_{mj}$ converted to $P$-values can be thought of as using conditional $P$-values at stage $m$ for a centered variable $j$, conditioned on the other remaining variables, assuming all nulls are true.

**6. Treatments vs. control, change point and successive correlation models.** The first two models of this section entail independent random samples from $(M+1)$ normal populations. Let $Z_{ij}$, $i=1,\dots,M+1$, $j=1,\dots,n$, be $N(\nu_i,\sigma^2)$. In the treatments vs. control model, the treatments correspond to $i=1,\dots,M$, while the control population corresponds to the $(M+1)$st population. Let $X_i = \bar{Z}_i - \bar{Z}_{M+1}$, $i=1,\dots,M$, so that $\mathbf{X}$ is distributed as multivariate normal with mean vector $\boldsymbol{\mu}$, $\mu_i = \nu_i - \nu_{M+1}$ and covariance matrix $(2\sigma^2/n)\Sigma$, where $\Sigma$ is intraclass with diagonal elements 1 and off diagonal elements 1/2. Should $\sigma^2$ be unknown, then an unbiased estimator of $\sigma^2$ is

$$\begin{aligned}
s^2 &= \sum_{i=1}^{M+1}\sum_{j=1}^{n}(Z_{ij}-\bar{Z}_i)^2/(M+1)(n-1) \\
&= \left(\sum_{i=1}^{M+1}\sum_{j=1}^{n}Z_{ij}^2 - n\sum_{i=1}^{M+1}\bar{Z}_i^2\right)\Big/(M+1)(n-1) \\
&= \left(T - n\sum_{i=1}^{M+1}\bar{Z}_i^2\right)\Big/(M+1)(n-1).
\end{aligned}$$

Furthermore, $(M+1)(n-1)s^2/\sigma^2$ is distributed as chi-square with $(M+1)(n-1)$ degrees of freedom and is independent of $\bar{\mathbf{Z}}$ and hence $\mathbf{X}$. We recognize that, in terms of $\mathbf{X}$ and $T$, we have a special case of the model of Section 2 and, in fact, we have one of the models of Section 4. For this problem then, MRD is an admissible procedure for two-sided as well as one-sided alternatives. The LRSD procedure is admissible for one-sided procedures and for two-sided procedures for $M=2$ and 3. For $M=4$, however, many counterexamples to admissibility exist. In the next section, we use simulations to evaluate MRD and compare it to the popular step-wise procedures that are based on $P$-values from marginal distributions.

The model for the change point problem also entails $M+1$ independent random samples from normal populations. Let $Z_{ij}$ be as in the previous setup, only this time we are interested in $M$ null hypotheses $H_i : \mu_i = \nu_{i+1} - \nu_i = 0$ vs. $K_i : \mu_i \neq 0$ for two-sided alternatives or $K_i^* : \mu_i > 0$ for one-sided alternatives, $i=1,\dots,M$. Let $X_i = \bar{Z}_{i+1} - \bar{Z}_i$, so that $\mathbf{X}$ is distributed as multivariate normal with mean vector $\boldsymbol{\mu}$ and covariance matrix $(\sigma^2/n)\Sigma$, where $\Sigma = (\sigma_{ij})$, and

$$\sigma_{ii} = 2, \sigma_{ij} = -1, \qquad \text{if } |i-j| = 1,\ \sigma_{ij} = 0, \quad \text{otherwise},\ i,j = 1,\dots,M, i \neq j. \tag{6.1}$$

Note that a rejected $H_i$ amounts to infering that a change in mean has occurred from time $i$ to time $(i+1)$. One seeks to identify all change points. There is a substantial literature on the change point problem [see, e.g., Chen and Gupta (2000), where reference is made to the binary segmentation procedure (BSP) due to Vostrikova (1981)].

For this problem, one can consider a number of approaches. Among them are MRD, LRSD and BSP, the usual step-up and step-down procedures based on $P$-values. There is a very interesting connection between the MRD and BSP methods. Both are based on the $U_{mi}$ statistics given in (2.4). This is further support for our general methodology since, in this special case, our statistics are precisely the statistics $V_{mi}$ used by Vostrikova (1981) for the change point problem. MRD and BSP use the statistics differently.

We now demonstrate that $U_{mi}$ are the same as $V_{mi}$ and note that the $U_{mi}$ statistics can be computed readily for any size problem, since it will not be necessary to actually invert any matrix or submatrix of $\Sigma$ as given in (6.1). Toward this end, for $1 \leq p \leq M$, define the $p \times p$ matrix

$$\Sigma(p) = \begin{pmatrix} 2 & -1 & 0 & 0 & \cdots & 0 & 0 & 0 \\ -1 & 2 & -1 & 0 & \cdots & 0 & 0 & 0 \\ \vdots & \vdots & \vdots & \vdots & & \vdots & \vdots & \vdots \\ 0 & 0 & 0 & 0 & \cdots & -1 & 2 & -1 \\ 0 & 0 & 0 & 0 & \cdots & 0 & -1 & 2 \end{pmatrix}. \tag{6.2}$$

Note that $\Sigma(M)$ was given in (6.1). It is easily verified that the first and last rows of $\Sigma^{-1}(p)$ are $(1/(p+1))(p(p-1)\cdots 1)$ and $(1/(p+1))(12\cdots p)$, respectively.

Suppose now that, at stage $m$, we have not eliminated $x_1$ (i.e., not rejected $x_1$ at an earlier stage) but have eliminated $x_{j_1},\dots,x_{j_{m-1}}$ (i.e., rejected $H_{j_1},\dots,H_{j_{m-1}}$). For this development, we may take $j_1<j_2<\cdots<j_{m-1}$ without loss of generality. Let $r$ be an index, $0\le r\le m$, and let $j_0=1$, $j_m=M$. We are now ready to prove the following theorem.

THEOREM 6.1. *For $j_r<i<j_{r+1}$, the statistics $U_{mi}$ of (2.4) can be expressed as*

$$[(j_{r+1}-j_r)/(j_{r+1}-i)(i-j_r)]^{1/2} \times\left[\sum_{j=j_r+1}^{i} Z_j-(i-j_r)\left(\sum_{j=j_r+1}^{j_{r+1}} Z_j\right)\Big/(j_{r+1}-j_r)\right]. \tag{6.3}$$

PROOF. See Appendix. □

REMARK 6.1. It can be shown that the BSP procedure is also admissible.

REMARK 6.2. For the change point model when $M\ge 4$, it can be demonstrated that the LRSD method is frequently inadmissible both for two-sided and one-sided alternatives.

The successive correlation model starts with an $M\times 1$ random vector $\mathbf{X}$, which is multivariate normal with mean vector $\boldsymbol{\mu}$ and covariance matrix

$$\Sigma(M)=\begin{pmatrix} 1 & \rho & 0 & 0 & \cdots & 0 & 0\\ \rho & 1 & \rho & 0 & \cdots & 0 & 0\\ 0 & \rho & 1 & \rho & \cdots & 0 & 0\\ \vdots & \vdots & \vdots & \vdots & & \vdots & \vdots\\ 0 & 0 & 0 & 0 & \cdots & 1 & \rho\\ 0 & 0 & 0 & 0 & \cdots & \rho & 1 \end{pmatrix}. \tag{6.4}$$

Note that if $\Sigma(0)\equiv 1$, then for $r=0,1,\dots,M$,

$$|\Sigma(r)|=|\Sigma(r-1)|-\rho^2|\Sigma(r-2)|. \tag{6.5}$$

Also, one can verify that the first row of the inverse of $\Sigma(r)$ is $(d_{1,r},\dots,d_{r,r})$ where

$$d_{i,r}=(-\rho)^{i-1}|\Sigma(r-i)|/|\Sigma(r)|. \tag{6.6}$$

By symmetry, the last row is $(d_{r,r},\dots,d_{1r})$. Proceeding as in the change point model case when $x_{j_1},\dots,x_{j_{m-1}}$ have been eliminated and $j_r<i<j_{r+1}$, the numerator of $U_{mi}$ [see (A.7)] is

$$\begin{aligned} & x_i-(0,\dots,0,\rho,\rho,0,\dots,0) \\ & \quad\times\begin{pmatrix}\Sigma(j_1-1) & & & & & \\ & \ddots & & & & \\ & & \Sigma(s_1) & & & \\ & & & \Sigma(s_2) & & \\ & & & & \ddots & \\ & & & & & \Sigma(M-j_{m-1})\end{pmatrix} \\ & \quad\times\mathbf{x}^{(i,j_1,\dots,j_{m-1})}, \end{aligned}\tag{6.7}$$

where the row vector above is of order $(M-m)\times 1$ and has two entries of $\rho$ in positions $i-1$ and $i$. If $i=1$ or $M$, then there is only one entry of $\rho$. Defining $s_1$ and $s_2$ as in the change point model, we find that (6.7) is

$$x_i+(0,\dots,0,d_{i,s_1},\dots,d_{s_1,s_1},d_{s_2,s_2},\dots,d_{1,s_2},0,\dots,0)\mathbf{x}^{(i,j_1,\dots,j_{m-1})},\tag{6.8}$$

where the nonzero entries in (6.8) appear in positions $j_r+1, j_{r+1}-2$. Thus, (6.8) becomes

$$x_i+\sum_{j=1}^{s_1}d_{j,s_1}x_{j_r+j}+\sum_{j=1}^{s_2}d_{s_2-j+1,s_2}x_{i+j}.\tag{6.9}$$

The denominator of $U_{mi}$ is

$$\begin{aligned} & 1-(0,\dots,0,\rho,\rho,0,\dots,0) \\ & \quad\times\begin{pmatrix}\Sigma(j_1-1) & & & & & \\ & \ddots & & & & \\ & & \Sigma(s_1) & & & \\ & & & \Sigma(s_2) & & \\ & & & & \ddots & \\ & & & & & \Sigma(M-j_m-1)\end{pmatrix}\begin{pmatrix}0\\ \vdots\\ 0\\ \rho\\ \rho\\ 0\\ \vdots\\ 0\end{pmatrix}, \end{aligned}\tag{6.10}$$

where the vectors are $(M-m)$-dimensional with two entries of $\rho$ in the $(i-1)$ and $i$ positions. If $i=1$ or $M$, then there is only one entry of $\rho$. Thus, only the $(s_1,s_1)$ element of $\Sigma^{-1}(s_1)$ and the $(1,1)$ element of $\Sigma^{-1}(s_2)$ will be needed. Specifically, we get

$$1-\rho^2\biggl(\frac{|\Sigma(s_1-1)|}{|\Sigma(s_1)|}+\frac{|\Sigma(s_2-1)|}{|\Sigma(s_2)|}\biggr).\tag{6.11}$$

**7. Simulations.** The MRD procedure can be viewed as a family of admissible procedures parametrized by a set of constants $\{C_1, \ldots, C_M\}$. It can be shown, using an inequality due to Sidák (1968), that $\{C_1, \ldots, C_M\}$ can be chosen so that the MRD procedure controls the strong FWER. However, such a choice of $C$'s would be extremely conservative and would sacrifice the gains achieved by MRD, which takes advantage of the correlation among the variables. It is also possible to choose smaller $C$'s to control FDR. However, this too is likely to lead to an overly conservative procedure. To determine a reasonable set of constants, one must study the risks (errors and error rates) for various choices of constants. As is the case in a typical decision theory problem where no optimal procedure exists, one must choose from a number of admissible procedures. This process needs to be done prior to looking at the data. To make this choice in practice, one must consider the particular application. In the examples we present, the number of hypotheses is very large, and so one expects only a small or modest percentage of alternatives to be true. Thus, we focused on that portion of the parameter space where, at most, 25% of alternatives were true. A large variety of sets of constants were evaluated through simulation. Those presented gave a good balance of performance in terms of expected numbers of type I and type II errors committed.

We have seen, in Section 3, that the MRD procedures possess the intuitive convexity property needed for admissibility regardless of the covariance matrix, $\sigma^2\Sigma$. These stepwise procedures make extensive use of the covariance structure at every stage. To see the types of improvements that can be made over usual stepwise methods, we now present some simulation studies. In this section, we report results for the treatments versus control model in both the $\sigma^2$ known and unknown cases. We also look at the change point model.

Our studies focus on the situation where the number of populations is large and the number of true alternatives is less than 25%. For two-sided alternatives, we present a comparison of the MRD method with either the step-up or step-down method (whichever did best in the given situation). The step-up and step-down methods used in the comparison are those based on $P$-values determined from marginal distributions. We report the expected number of type I errors, the expected number of type II errors and the FDR. To obtain the probabilities of types I and II errors we can divide the expected number of errors, in the tables below by the number of true nulls and alternatives, respectively. For all simulations, we used 1000 iterations.

Table 1 gives the results for the treatments versus control model (i.e., intraclass with a correlation coefficient of 0.5) with a known $\sigma^2 = 1$ for two-sided alternatives. The step-up procedure in the table is the Benjamini–Hochberg (1995) FDR controlling procedure where $\mathrm{FDR} = 0.05$. Thus, the critical values for the step-up procedure are $[\Phi^{-1}(0.05i/2M)]$. The critical values for MRD are somewhat related to the FWER controlling step-down procedure where the control is at level 0.05. Specifically, these critical values for MRD are as follows. For $\alpha = 0.05$, $M = 10{,}000$, $C_1 = \Phi^{-1}(1-\alpha/2M)$, $C_i = 0.71\Phi^{-1}(1-\alpha/2(M-i+1))$ and $1 < i \le M$. These critical values were selected by trial and error using simulations with 1000 iterations. They were chosen so that a desirable procedure would ensue and also to suggest a way to get critical values in other cases. Another consideration was to try to match step-up in FDR when the number of nonnulls is large. Here, $M = 10{,}000$, and the results are most dramatic. There is improvement (usually substantial) in both the expected number of types I and II errors.

Table 2 also gives results for the treatments versus control model (i.e., intraclass with a correlation coefficient of 0.5) but with unknown $\sigma^2$ for two-sided alternatives. Here, $M = 3000$, $n = 10$, $\alpha = 0.05$, $C_1 = \Phi^{-1}(1-\alpha/2M)$, $C_i = 0.63\Phi^{-1}(1-\alpha/2(M-i+1))$ and $1 < i \le M$. The step-down procedure in the table is based on $P$-values of the marginal distributions of $t$-statistics with a pooled estimate of $\sigma^2$. The step-down procedure controls FWER at $\alpha = 0.05$.

Table 3 deals with the change point model for two-sided alternatives. Unlike the intraclass model, the variables are not exchangeable. Thus, the pattern of true mean values as well as the choice of true mean values impacts the operating characteristics of the procedures. It would be difficult to select a particular portion of the parameter space to study without knowing the specific application. The type of pattern in mean values we present reflects the notion of an occasional rise in mean value as follows. The sequence of differences in consecutive means are of the form $0,\dots,0,1,1,1,0,\dots,0,1,1,1,0,\dots,0$ where the sets of triples $(1,1,1)$ are equally spaced. Once again, $M = 3000$, $\alpha = 0.05$, $C_1 = \Phi^{-1}(1-\alpha/2M)$, $C_i = 0.77\Phi^{-1}(1-\alpha/2(M-i+1))$. The step-down procedure in the table is based on the difference of two normal variables, each with variance 1. The procedure controls FWER at $\alpha = 0.05$.

The message in Tables 2 and 3 for two sided alternatives is that MRD has slightly higher expected number of type I errors but has many fewer type II errors.

Table 4 gives results for the treatments versus control model with known $\sigma^2 = 1$, for one-sided alternatives. We compare MRD, LRSD, step-down based on Dunnett's tests for $\alpha = 0.05$, call it $D(0.05)$, step-down based on Dunnett's tests for $\alpha = 0.2$, call it $D(0.2)$, regular step-down (SD) and regular step-up (SU). Before commenting on the simulation findings, we make some remarks. MRD and LRSD both take dependency into account in two ways. Namely, through test statistics and through critical values. $D(0.05)$ and $D(0.2)$ take dependency into account only through critical values, SD and SU do not take dependency into account at all. Recall that our proposal is to sacrifice some FDR control, especially when there are not too

Table 1

*Comparison of MRD and SU procedures for treatments vs. control, variance known*

| **Number of means equal to** | | | | | **Expected # of type I errors** | | **Expected # of type II errors** | | **FDR** | | **Total errors** | |
|---|---|---|---|---|---|---|---|---|---|---|---|---|
| **0** | **−4** | **−2** | **2** | **4** | **MRD** | **SU** | **MRD** | **SU** | **MRD** | **SU** | **MRD** | **SU** |
| 10000 | 0 | 0 | 0 | 0 | 0.67 | 28 | 0 | 0 | 0.05 | 0.02 | 0.67 | 28 |
| 9200 | 0 | 800 | 0 | 0 | 13.02 | 24.03 | 560.32 | 726.5 | 0.05 | 0.02 | 573.34 | 750.52 |
| 9200 | 800 | 0 | 0 | 0 | 12.23 | 58.77 | 5.99 | 131.18 | 0.02 | 0.04 | 18.22 | 189.96 |
| 8400 | 0 | 800 | 800 | 0 | 11.2 | 40.32 | 1041.91 | 1463.22 | 0.02 | 0.03 | 1053.11 | 1503.54 |
| 8400 | 0 | 0 | 1600 | 0 | 16.06 | 43.45 | 1205.59 | 1392.09 | 0.04 | 0.02 | 1221.65 | 1435.54 |
| 8400 | 800 | 0 | 800 | 0 | 12.78 | 55.09 | 557.82 | 730.51 | 0.01 | 0.03 | 570.60 | 785.6 |
| 8400 | 0 | 0 | 800 | 800 | 12.95 | 34.40 | 563.96 | 752.64 | 0.01 | 0.03 | 576.91 | 787.04 |
| 8400 | 800 | 0 | 0 | 800 | 13.28 | 73.65 | 12.37 | 148.81 | 0.01 | 0.04 | 25.65 | 222.45 |
| 8400 | 0 | 0 | 0 | 1600 | 13.46 | 70.82 | 12.56 | 167.88 | 0.01 | 0.04 | 26.02 | 238.7 |
| 7600 | 0 | 800 | 1600 | 0 | 12.17 | 55.13 | 1602.92 | 2121.25 | 0.02 | 0.03 | 1614.47 | 2176.37 |
| 7600 | 0 | 0 | 2400 | 0 | 24.95 | 59.77 | 1943.43 | 2000.7 | 0.05 | 0.03 | 1968.37 | 2060.47 |
| 7600 | 800 | 0 | 1600 | 0 | 16.17 | 57.67 | 1191.87 | 1313.02 | 0.01 | 0.03 | 1208.05 | 1370.7 |
| 7600 | 0 | 0 | 1600 | 800 | 16.41 | 58.33 | 1202.26 | 1326.52 | 0.01 | 0.03 | 1218.66 | 1384.85 |
| 7600 | 800 | 0 | 800 | 800 | 14.32 | 85.26 | 562.51 | 718.44 | 0.01 | 0.03 | 576.83 | 803.7 |
| 7600 | 0 | 0 | 800 | 1600 | 14.56 | 69.92 | 569.23 | 758.13 | 0.01 | 0.04 | 583.8 | 828.05 |
| 7600 | 800 | 0 | 0 | 1600 | 14.73 | 95.19 | 19.79 | 160.22 | 0.01 | 0.03 | 34.52 | 255.4 |
| 7600 | 0 | 0 | 0 | 2400 | 15.58 | 116.56 | 21.17 | 218.25 | 0.01 | 0.04 | 36.76 | 334.82 |

TABLE 2

*Comparison of MRD and SD procedures for treatments vs. control, variance unknown*

| **Number of noncentrality parameters equal to the value** | | | | | **Expected # of type I errors** | | **Expected # of type II errors** | | **FDR** | | **Total errors** | |
|---|---|---|---|---|---|---|---|---|---|---|---|---|
| **0** | **−3** | **−1** | **1** | **3** | **MRD** | **SD** | **MRD** | **SD** | **MRD** | **SD** | **MRD** | **SD** |
| 3000 | 0 | 0 | 0 | 0 | 0.18 | 0.09 | 0 | 0 | 0.02 | 0.02 | 0.18 | 0.09 |
| 2800 | 0 | 200 | 0 | 0 | 1.22 | 0.04 | 198.62 | 199.93 | 0.05 | 0.02 | 199.85 | 199.96 |
| 2800 | 200 | 0 | 0 | 0 | 7.33 | 0.07 | 22.7 | 180.9 | 0.04 | 0.01 | 30.03 | 180.97 |
| 2600 | 0 | 200 | 200 | 0 | 2.51 | 0.04 | 992.62 | 399.78 | 0.07 | 0.01 | 395.13 | 399.82 |
| 2600 | 0 | 0 | 400 | 0 | 1.7 | 0.04 | 396.81 | 399.79 | 0.05 | 0.01 | 398.51 | 399.83 |
| 2600 | 200 | 0 | 200 | 0 | 6.5 | 0.03 | 211.3 | 379.14 | 0.03 | 0 | 217.81 | 379.17 |
| 2600 | 0 | 0 | 200 | 200 | 7.02 | 0.04 | 218.62 | 381.31 | 0.04 | 0.01 | 225.64 | 381.35 |
| 2600 | 200 | 0 | 0 | 200 | 4.31 | 0.07 | 58.43 | 361.18 | 0.01 | 0 | 62.75 | 361.25 |
| 2600 | 0 | 0 | 0 | 400 | 4.92 | 0.04 | 59.52 | 362.74 | 0.01 | 0.01 | 64.43 | 362.78 |
| 2400 | 0 | 200 | 400 | 0 | 2.82 | 0.05 | 587.53 | 599.69 | 0.06 | 0.01 | 590.35 | 599.74 |
| 2400 | 0 | 0 | 600 | 0 | 1.74 | 0.05 | 596.39 | 599.66 | 0.05 | 0.01 | 598.14 | 599.71 |
| 2400 | 200 | 0 | 400 | 0 | 6.2 | 0.02 | 403.79 | 580.5 | 0.03 | 0 | 409.99 | 580.52 |
| 2400 | 0 | 0 | 400 | 200 | 7.07 | 0.02 | 417.96 | 581.49 | 0.04 | 0.01 | 425.03 | 581.51 |
| 2400 | 200 | 0 | 200 | 200 | 3.88 | 0.03 | 254.02 | 562.3 | 0.01 | 0 | 257.9 | 562.32 |
| 2400 | 0 | 0 | 200 | 400 | 5.01 | 0.02 | 262.95 | 562.79 | 0.02 | 0.01 | 267.96 | 562.81 |
| 2400 | 200 | 0 | 0 | 400 | 2.91 | 0.05 | 110.07 | 541.56 | 0.01 | 0 | 112.98 | 541.61 |
| 2400 | 0 | 0 | 0 | 600 | 3.96 | 0.03 | 110.39 | 543.14 | 0.01 | 0.01 | 114.36 | 543.17 |

Table 3

*Comparison of MRD and SD procedures for the change point model*

| Number of | | Expected # of type I errors | | Expected # of type II errors | | FDR | | Total errors | |
|---|---|---|---|---|---|---|---|---|---|
| nulls | triples | MRD | SD | MRD | SD | MRD | SD | MRD | SD |
| 3000 | 0 | 0 | 0.05 | 0 | 0 | 0 | 0.05 | 0 | 0.05 |
| 2970 | 10 | 1.81 | 0.04 | 21.14 | 30 | 0.16 | 0.04 | 22.95 | 30.04 |
| 2955 | 15 | 2.04 | 0.06 | 31.34 | 44.99 | 0.13 | 0.06 | 33.38 | 45.05 |
| 2925 | 25 | 4.27 | 0.05 | 52.65 | 74.98 | 0.16 | 0.05 | 56.93 | 75.03 |
| 2850 | 50 | 8.39 | 0.05 | 105.26 | 149.97 | 0.16 | 0.05 | 113.65 | 150.03 |
| 2820 | 60 | 7.63 | 0.06 | 125.77 | 179.97 | 0.12 | 0.05 | 132.99 | 180.03 |
| 2700 | 100 | 17.61 | 0.04 | 210.99 | 299.95 | 0.17 | 0.04 | 228.6 | 299.99 |
| 2550 | 150 | 27.21 | 0.05 | 317.71 | 449.93 | 0.17 | 0.04 | 344.92 | 449.97 |
| 2400 | 200 | 36.52 | 0.04 | 423.90 | 599.90 | 0.17 | 0.04 | 460.42 | 599.95 |

many rejections. MRD was recommended when the proportion of false nulls is less than 0.25. Also, recall the geometric relationship between MRD and LRSD, as noted in Section 5. In light of this, we expect and do observe that the performance of MRD and LRSD (in terms of expected number of mistakes) is comparable. One advantage that MRD has over LRSD is in computation. LRSD requires a package like "quadprog" in R. This program is very time consuming for $M > 100$, which is why the simulation is done for $M = 100$ and not for a larger $M$. Since $D(0.05)$ takes dependency into account through critical values, that procedure should and does perform better than SD, which controls FWER at $\alpha = 0.05$. It is not fair in a sense to compare $D(0.05)$, a markedly conservative procedure, with MRD and LRSD. However, one can compare $D(0.2)$ with MRD and LRSD, and the latter two are preferred.

The simulations are based on 1000 iterations. The largest percentage of true alternatives considered is 25. The cirtical values for MRD and LRSD are as follows. Let $C_1(SD) = \Phi^{-1}(1 - \alpha/M)$, $C_i(SD) = \Phi^{-1}(1 - \alpha/(M - i + 1))$,$1 < i \leq M$, then $C_1$ for LRSD is $1.25C_1(SD)$ and $C_i$ for LRSD is $1.2C_i(SD)$, for $i \neq 1$. For MRD $C_1 = C_1(SD)$, $C_i = 0.7C_i(SD)$, for $i \neq 1$. The $C$'s for $D(0.2)$ and $D(0.05)$ are obtained by simulation.

Table 4 offers simulations that yield FDR and total errors for each of the six procedures. Other simulations yielded expected number of type I errors and expected number of type II errors. These are not given in the table, because the pattern for type I errors is the same as with FDR, and the pattern for type II errors can be discerned from the columns giving total errors.

Table 4

*Comparison of MRD, LRSD, $D(0.2)$, $D(0.05)$, SD and SU for one-sided treatments vs. control*

| # of means equal to | | | FDR | | | | | | Total errors | | | | | |
|---|---|---|---|---|---|---|---|---|---|---|---|---|---|---|
| **0** | **2** | **4** | **MRD** | **LRSD** | $\boldsymbol{D(0.2)}$ | $\boldsymbol{D(0.05)}$ | **SD** | **SU** | **MRD** | **LRSD** | $\boldsymbol{D(0.2)}$ | $\boldsymbol{D(0.05)}$ | **SD** | **SU** |
| 100 | 0 | 0 | 0.05 | 0.04 | 0.19 | 0.05 | 0.03 | 0.04 | 0.11 | 0.19 | 0.95 | 0.08 | 0.07 | 0.63 |
| 95 | 5 | 0 | 0.1 | 0.11 | 0.11 | 0.02 | 0.01 | 0.02 | 4.53 | 4.81 | 4.99 | 4.89 | 4.91 | 5.19 |
| 95 | 0 | 5 | 0.17 | 0.06 | 0.08 | 0.02 | 0.01 | 0.04 | 1.39 | 1.54 | 2.46 | 3.11 | 3.33 | 3.40 |
| 90 | 5 | 5 | 0.12 | 0.07 | 0.05 | 0.01 | 0.01 | 0.03 | 4.40 | 5.13 | 6.77 | 7.89 | 8.38 | 8.23 |
| 90 | 10 | 0 | 0.09 | 0.14 | 0.09 | 0.02 | 0.01 | 0.04 | 8.28 | 8.12 | 9.18 | 9.65 | 9.73 | 10.06 |
| 90 | 0 | 10 | 0.1 | 0.04 | 0.05 | 0.01 | 0.00 | 0.04 | 1.66 | 2.20 | 4.29 | 6.09 | 6.82 | 5.93 |
| 85 | 5 | 10 | 0.07 | 0.04 | 0.04 | 0.01 | 0.00 | 0.04 | 4.53 | 5.51 | 8.24 | 10.58 | 11.51 | 10.02 |
| 85 | 10 | 5 | 0.08 | 0.07 | 0.05 | 0.01 | 0.01 | 0.04 | 7.56 | 8.24 | 10.96 | 12.55 | 13.11 | 12.70 |
| 80 | 15 | 5 | 0.06 | 0.06 | 0.04 | 0.01 | 0.00 | 0.04 | 10.86 | 11.21 | 14.81 | 17.34 | 17.86 | 16.75 |
| 80 | 5 | 15 | 0.05 | 0.02 | 0.02 | 0.01 | 0.00 | 0.03 | 4.78 | 6.03 | 10.00 | 13.27 | 15.14 | 11.82 |
| 80 | 10 | 10 | 0.06 | 0.04 | 0.03 | 0.00 | 0.00 | 0.03 | 7.85 | 8.66 | 12.7 | 15.50 | 16.55 | 14.63 |
| 80 | 20 | 0 | 0.04 | 0.10 | 0.05 | 0.01 | 0.01 | 0.03 | 15.88 | 13.54 | 17.31 | 19.09 | 19.34 | 19.00 |
| 80 | 0 | 20 | 0.05 | 0.01 | 0.02 | 0.00 | 0.00 | 0.03 | 2.06 | 3.34 | 7.50 | 11.66 | 13.76 | 9.09 |
| 75 | 5 | 20 | 0.04 | 0.02 | 0.02 | 0.00 | 0.00 | 0.03 | 5.23 | 6.68 | 11.49 | 16.59 | 18.43 | 13.06 |
| 75 | 20 | 5 | 0.04 | 0.06 | 0.03 | 0.01 | 0.00 | 0.03 | 14.51 | 14.17 | 18.95 | 22.10 | 22.73 | 21.06 |
| 75 | 15 | 10 | 0.05 | 0.04 | 0.03 | 0.01 | 0.00 | 0.03 | 11.02 | 11.66 | 16.55 | 20.21 | 21.34 | 18.62 |
| 75 | 10 | 15 | 0.05 | 0.03 | 0.02 | 0.01 | 0.00 | 0.03 | 8.13 | 9.26 | 14.01 | 18.26 | 19.82 | 15.81 |
| 75 | 25 | 0 | 0.03 | 0.08 | 0.04 | 0.01 | 0.00 | 0.03 | 20.35 | 16.52 | 21.73 | 23.87 | 24.28 | 23.63 |
| 75 | 0 | 25 | 0.05 | 0.01 | 0.02 | 0.00 | 0.00 | 0.03 | 2.50 | 3.88 | 9.14 | 14.34 | 17.12 | 10.59 |

In summary, under the stated conditions, the admissible procedures MRD and LRSD have comparable performances with a computational advantage for MRD, should $M > 100$.

## APPENDIX

PROOF OF LEMMA 3.2. For $j \neq 1$, $j \neq j_1, \dots, j \neq j_{(m-1)}$, use (2.4) and recall $\mathbf{g}$ is the first column of $\Sigma$ to see that

$$\begin{aligned}
&U_{mj}(\mathbf{x} + r\mathbf{g}) \\
&\quad = \{x_j + r\sigma_{j1} - \boldsymbol{\sigma}_j^{(j_1,\dots,j_{(m-1)})'} \\
&\quad\quad \times \Sigma^{-1}_{(j_1,\dots,j_{(m-1)},j)}(\mathbf{x}^{(j_1,\dots,j_{(m-1)},j)} + r\mathbf{g}^{(j_1,\dots,j_{(m-1)},j)})\}/\sigma_{(j\cdot j_1,\dots,j_{(m-1)})} \qquad \text{(A.1)} \\
&\quad = U_{mj}(\mathbf{x}) + [r\sigma_{j1} - r\boldsymbol{\sigma}_j^{(j_1,\dots,j_{(m-1)})'}(1, 0, \dots, 0)']/\sigma_{(j\cdot j_1,\dots,j_{(m-1)})} \\
&\quad = U_{mj}(\mathbf{x}).
\end{aligned}$$

This establishes (3.4).

Now,

$$\begin{aligned}
&U_{m1}(\mathbf{x} + r\mathbf{g}) \\
&\quad = U_{m1}(\mathbf{x}) + [r\sigma_{11} - r\sigma_{(1)}^{(j_1,\dots,j_{(m-1)})'}\Sigma^{-1}_{(j_1,\dots,j_{(m-1)},1)}\sigma_{(1)}^{(j_1,\dots,j_{(m-1)})}] \qquad \text{(A.2)} \\
&\quad\quad /\sigma_{(1\cdot j_1,\dots,j_{(m-1)})} \\
&\quad = U_{m1}(\mathbf{x}) + r,
\end{aligned}$$

which establishes (3.5). $\square$

PROOF OF LEMMA 3.3. If $\phi_U(\mathbf{x}^*) = 0$ then, when $\mathbf{x}^*$ is observed, the process must stop before $H_1$ is rejected. Suppose it stops at stage $m$ without having rejected $H_1$. That means that $U_{m,j_m} < C_m$, which is equivalent to $|U_{mi}| < C_m$ for all $i \in \{1, \dots, M\} \setminus j_1, \dots, j_{m-1}$, $j_i \neq 1$. Also $U_{ij_1} \geq C_i$, $i = 1, \dots, m-1$, $j_i \neq 1$. Next, consider $\mathbf{x}^* + r_0\mathbf{g}$, which is a reject $H_1$ point. By Lemma 3.2, (A.1) and (A.2) imply that only the functions $U_{m1}(\mathbf{x})$ can change from $\mathbf{x}^*$ to $\mathbf{x}^* + r_0\mathbf{g}$. Also, at some stage $h \leq m$ from (A.2), $U_{h,1}$ must have increased to become positive and become the maximum function at that stage and also be $\geq C_h$. By (A.2), $U_{h,1}$ will be at least this large for all $r \geq r_0$. Thus, $H_1$ will also be rejected for all $\mathbf{x} + r\mathbf{g}$, $r > r_0$. $\square$

PROOF OF THEOREM 4.1. We prove the theorem for $M = 3$. For $M = 2$ the method is the same and the proof is simpler. In light of Lemma 3.1 and the proof of Theorem 3.1, we need to show that the LRSD test for $H_1$ vs. $K_1$, say $\phi_1(\mathbf{x})$, as a function of $\mathbf{x}+\gamma\mathbf{g}$, goes from reject to accept to reject as $\gamma$ varies from $(-\infty,\infty)$. Another way of stating this requirement is, suppose $\phi_1(\mathbf{x}^*)=1$ when $x_1^*>0$. Then, we must have $\phi_1(\mathbf{x}^*+\gamma\mathbf{g})=1$ for all $\gamma>0$, while if $\phi_1(\mathbf{x}^*)=1$ and $x_1^*<0$, we must have $\phi_1(\mathbf{x}^*-\gamma\mathbf{g})=1$ for all $\gamma>0$.

There are a number of cases that need to be treated. Namely each of the three stages at which $H_1$ is rejected. If $H_1$ is rejected at stage 1 at $\mathbf{x}=\mathbf{x}^*$ and $x_1^*>0$, with $x_1^*>|x_2^*|$ and $x_1^*>|x_3^*|$, then $\mathbf{x}^{*\prime}\Sigma^{-1}\mathbf{x}^*\geq C_1$, and this implies that

$$(\mathbf{x}^*+\gamma\mathbf{g})'\Sigma^{-1}(\mathbf{x}^*+\gamma\mathbf{g})=\mathbf{x}^{*\prime}\Sigma^{-1}\mathbf{x}^*+2\gamma x_1^*+\gamma^2>C_1.$$

Also, $x_1^*+\gamma>|x_2^*+\gamma\rho|$ and $x_1^*+\gamma>|x_3^*+\gamma\rho|$, which means that $\phi_1(\mathbf{x}^*)=1$ implies $\phi_1(\mathbf{x}^*+\gamma\mathbf{g})=1$ and $H_1$ is rejected at stage 1 for all $\mathbf{x}^*+\gamma\mathbf{g}$, all $\gamma>0$.

The next case to consider is when $H_1$ is rejected at the second stage for $\mathbf{x}=\mathbf{x}^*$. Two subcases are $x_1^*>0$ and $x_1^*<0$. For $x_1^*>0$, suppose $x_3$ is out first. Then, we find that ${x_1^*}^2+{x_2^*}^2-2\rho x_1^*x_2^*\geq C_2$ and

$$\begin{aligned}(A.3)\qquad &(x_1^*+\gamma)^2+(x_2^*+\gamma\rho)^2-2\rho(x_1^*+\gamma)(x_2^*+\rho\gamma)\\ &\quad={x_1^*}^2+{x_2^*}^2-2\rho x_1^*x_2^*+2\gamma x_1^*+\gamma^2+\rho^2\gamma^2-2\rho^2x_1^*\gamma-2\rho^2\gamma^2.\end{aligned}$$

But, since $\gamma^2+\rho^2\gamma^2>2\rho^2\gamma^2$ and $2\gamma x_1^*\geq 2\rho^2\gamma x_1^*$, it follows that (A.3) $>C_2$ for all $\gamma>0$. Hence, $\phi_1(\mathbf{x}^*+\gamma\mathbf{g})=1$ for all $\gamma>0$. If $x_1^*<0$, a similar argument works for $(x_1^*-\gamma)^2+(x_2^*-\gamma\rho)^2-2\rho(x_1^*-\gamma)(x_2^*-\gamma\rho)$.

Finally, the third case is when $H_1$ is rejected at stage 3. In subcases where the ordering of the components of $\mathbf{x}^*$ is maintained with $(\mathbf{x}^*+\gamma\mathbf{g})$, it is easy to prove the required monotonicity property. The most challenging subcase is if $|x_3^*|>x_2^*>x_1^*>0$ with $x_3^*<0$ but

$$(A.4)\qquad |x_3^*+\gamma\rho|<x_2^*+\gamma\rho.$$

In this case, when $\rho>0$, we use the fact that ${x_3^*}^2>{x_2^*}^2$ and use inequalities, as in the previous case, to prove the result. When $\rho<0$, we observe that, if $|x_3^*|>x_2^*>x_1^*>0$ and $x_3^*<0$, then $|x_3^*+\rho\gamma|>|x_2^*+\rho\gamma|$, and so (A.4) cannot happen. It is easy to verify then that if $\phi_1(\mathbf{x}^*)=1$ then $\phi_1(\mathbf{x}^*+\gamma\mathbf{g})=1$ for all $\gamma>0$. Similarly, for $x_1<0$. □

Proof of Theorem 4.2. Once again, we focus on $H_1$ vs. $K_1^*$ and demonstrate that if $\phi_1(\mathbf{x}^*)=1$, then $\phi_1(\mathbf{x}^*+\gamma\mathbf{g})=1$ for all $\gamma>0$. Suppose $H_1$ is rejected at stage $m$ at $\mathbf{x}=\mathbf{x}^*$. Then, $x_{j_1}^*>x_{j_2}^*>\cdots>x_{j_{m-1}}^*>x_1^*>x_{j_{m+1}}^*>\cdots>x_{j_M}^*$ and $x_1^*>0$. Note that, at $\mathbf{x}^{**}=\mathbf{x}^*+\gamma\mathbf{g}$, the orders of all coordinates are preserved except perhaps the first coordinate, which now can be anywhere among the $m$ largest coordinates. The $k$ stage global hypothesis is considered if $H_{j_1},\dots,H_{j_{k-1}}$ have been rejected. This global testing problem is $H_{kG}: \boldsymbol{\mu}^{(j_1,\dots,j_{k-1})} = \mathbf{0}^{(j_1,\dots,j_{k-1})}$ vs. $K_{kG}: \boldsymbol{\mu}^{(j_1,\dots,j_{k-1})} \in Q^{(j_1,\dots,j_{k-1})}$, where $Q^{(j_1,\dots,j_{k-1})} = \{\boldsymbol{\mu}^{(j_1,\dots,j_{k-1})}: \mu_i \geq 0, i \in [1,\dots,M]\setminus[j_1,\dots,j_{k-1}]\}\setminus \mathbf{0}^{(j_1,\dots,j_{k-1})}$. The likelihood ratio test rejects $H_{kG}$ if $\mathbf{x}^*$ is observed and if

$$\begin{aligned}
&\sup_Q \exp\{\mathbf{x}^{*(j_1,\dots,j_{k-1})'}\Sigma^{-1}_{(j_1,\dots,j_{k-1})}\boldsymbol{\mu}^{(j_1,\dots,j_{k-1})} \\
&\qquad - (1/2)\boldsymbol{\mu}^{(j_1,\dots,j_{k-1})'}\Sigma^{-1}_{(j_1,\dots,j_{k-1})}\boldsymbol{\mu}^{(j_1,\dots,j_{k-1})}\} \\
&\quad = \exp\{\mathbf{x}^{*(j_1,\dots,j_{k-1})'}\Sigma^{-1}_{(j_1,\dots,j_{k-1})}\hat{\boldsymbol{\mu}}^{(j_1,\dots,j_{k-1})*} \\
&\qquad - (1/2)\hat{\boldsymbol{\mu}}^{(j_1,\dots,j_{k-1})*'}\Sigma^{-1}_{(j_1,\dots,j_{k-1})}\hat{\boldsymbol{\mu}}^{(j_1,\dots,j_{k-1})*}\},
\end{aligned} \tag{A.5}$$

where $\hat{\boldsymbol{\mu}}^{(j_1,\dots,j_{k-1})*}$ is the maximum likelihood estimator of $\boldsymbol{\mu}^{(j_1,\dots,j_{k-1})}$, when $\mathbf{x} = \mathbf{x}^*$.

Next, consider the likelihood ratio test statistic at $\mathbf{x}^{**}$. It is

$$\begin{aligned}
&\sup_Q \exp\{(\mathbf{x}^{*(j_1,\dots,j_{k-1})} + \gamma\mathbf{g}^{(j_1,\dots,j_{k-1})})' \\
&\qquad \times \Sigma^{-1}_{(j_1,\dots,j_{k-1})}\boldsymbol{\mu}^{(j_1,\dots,j_{k-1})} \\
&\qquad - (1/2)\boldsymbol{\mu}^{(j_1,\dots,j_{k-1})'}\Sigma^{-1}_{(j_1,\dots,j_{k-1})}\boldsymbol{\mu}^{(j_1,\dots,j_{k-1})}\} \\
&\quad \geq \exp\{(\mathbf{x}^{*(j_1,\dots,j_{k-1})})' \\
&\qquad \times \Sigma^{-1}_{(j_1,\dots,j_{k-1})}\hat{\boldsymbol{\mu}}^{(j_1,\dots,j_{k-1})*} \\
&\qquad - (1/2)\hat{\boldsymbol{\mu}}^{(j_1,\dots,j_{k-1})*}\Sigma_{(j_1,\dots,j_{k-1})}\hat{\boldsymbol{\mu}}^{(j_1,\dots,j_{k-1})*}\}.
\end{aligned} \tag{A.6}$$

Recognize that the right-hand side of (A.6) is the maximized likelihood in (A.5) times $\exp\gamma\hat{\mu}_1^{(j_1,\dots,j_{k-1})*}$. Since $\hat{\mu}_1^{(j_1,\dots,j_{k-1})*} > 0$, it follows from (A.5) and (A.6) that (A.6) $\geq C_k$, which means that there is a rejection at stage $k$ at $\mathbf{x}^{**}$ if there was a rejection at stage $k$ at $\mathbf{x}^*$, $k = 1,\dots,M$. Since the order of the coordinates of $x_{j_1},\dots,x_{j_{m-1}}$ remains unchanged and $x_1^{**}$ is among the $m$ largest coordinates of $\mathbf{x}^{**}$, it follows that $H_1$ is rejected at stage $m$ or sooner. □

Proof of Theorem 6.1. With $x_{j_1},\dots,x_{j_{m-1}}$ eliminated, the covariance matrix of the remaining variables is the block matrix

$$\begin{pmatrix}
\Sigma(j_1-1) & & & 0 \\
& \Sigma(j_2-j_1-1) & & \\
& & \ddots & \\
0 & & & \Sigma(M-j_{m-1})
\end{pmatrix},$$

where $\Sigma(p)$ is given in (6.2). Now, let $s_1 = i - j_r - 1$ and $s_2 = j_{r+1} - i - 1$. Then, the residual for $x_i$ (i.e., the numerator of $U_{mi}$) is

$$x_i - (0, \ldots, 0, -1, -1, 0, \ldots, 0) \times \begin{pmatrix} \Sigma(j_1 - 1) & & & & & 0 \\ & \ddots & & & & \\ & & \Sigma(s_1) & & & \\ & & & \Sigma(s_2) & & \\ & & & & \ddots & \\ 0 & & & & & \Sigma(M - j_{m-1}) \end{pmatrix}^{-1} \times \mathbf{x}^{(i, j_1, \ldots, j_{m-1})}, \tag{A.7}$$

where the row vector above is of order $(M - m) \times 1$ and has two entries of $-1$ in positions $(i-1)$ and $i$. If $i = 1$ or $M$, then there is only one entry of $-1$. Thus, only the last row of $\Sigma^{-1}(s_1)$ and first row of $\Sigma^{-1}(s_2)$ will be needed. Specifically, (A.7) can be written as

$$\begin{aligned} &x_i + (0, \ldots, 0, 1/(s_1 + 1), 2/(s_1 + 1), \ldots, \\ &\qquad s_1/(s_1 + 1), s_2/(s_2 + 1), \ldots, 1/(s_2 + 1), 0, \ldots, 0)\mathbf{x}^{(i, j_1, \ldots, j_{m-1})}, \end{aligned} \tag{A.8}$$

where the nonzero entries in (A.8) appear in positions $j_r + 1, \ldots, j_{r+1} - 2$, so that the residual depends only on $x_{j_r+1}, \ldots, x_{j_{r+1}-1}$. Thus, (A.8) becomes

$$x_i + (1/(s_1 + 1)) \sum_{j=1}^{s_1} j x_{j_r + j} + (1/(s_2 + 1)) \sum_{j=1}^{s_2} (s_2 - j + 1) x_{i+j}. \tag{A.9}$$

Since $x_j = \bar{Z}_j - \bar{Z}_{j+1}$, (A.9) can be written as

$$\begin{aligned} &(1/(s_1 + 1)) \sum_{j=1}^{s_1+1} \bar{Z}_{j_r + j} - (1/(s_2 + 1)) \sum_{j=1}^{s_2+1} \bar{Z}_{i+j} \\ &\quad = (1/(i - j_r)) \sum_{j=j_r+1}^{i} \bar{Z}_j - (1/(j_{r+1} - i)) \sum_{j=i+1}^{j_{r+1}} \bar{Z}_j \\ &\quad = [(j_{r+1} - j_r)/(j_{r+1} - i)(i - j_r)] \\ &\qquad \times \left\{ \sum_{j=j_r+1}^{i} \bar{Z}_j - [(i - j_r)/(j_{r+1} - j_r)] \sum_{j=j_r+1}^{j_{r+1}} \bar{Z}_j \right\}. \end{aligned} \tag{A.10}$$

A similar computation yields the denominator of $U_{mi}$, namely $\sigma_{(i \cdot j_1, \ldots, j_{m-1})}$ as

$$[(j_{r+1} - j_r)/(j_{r+1} - i)(i - j_r)]^{1/2}. \tag{A.11}$$

Combine (A.10) and (A.11) and (6.3) is established. $\square$

A. Cohen
H. B. Sackrowitz
Department of Statistics
Rutgers University
Piscataway, New Jersey 08854
USA
E-mail: artcohen@rci.rutgers.edu
sackrowi@rci.rutgers.edu

M. Xu
Guanghua School of Management
Peking University
Beijing
China
E-mail: minyaxu@gsm.pku.edu.cn